\documentclass[11pt]{amsart}
\usepackage{epsfig}
\usepackage{amsfonts}
\usepackage{amsmath}
\usepackage{euscript}
\usepackage{amscd}
\usepackage[all]{xypic}

\def\and{\operatorname{and}}

\def\dim{\operatorname{dim}}
\def\depth{\operatorname{depth}}
\def\grade{\operatorname{grade}}

\def\ann{\operatorname{Ann}}

\def\modulo{\operatorname{modulo}}

\newcommand{\m}{\mathfrak m}

\newcommand{\R}{\mathcal R}

\newcommand{\F}{\mathcal F}

\newcommand{\bl}{\begin{lemma}}
\newcommand{\el}{\end{lemma}}
\newcommand{\bt}{\begin{theorem}}
\newcommand{\et}{\end{theorem}}
\newcommand{\ben}{\begin{enumerate}}
\newcommand{\een}{\end{enumerate}}
\newcommand{\bpf}{\begin{proof}}
\newcommand{\beqn}{\begin{eqnarray*}}
\newcommand{\eeqn}{\end{eqnarray*}}
\newcommand{\bd}{\begin{definition}}
\newcommand{\ed}{\end{definition}}
\newcommand{\bp}{\begin{proposition}}
\newcommand{\ep}{\end{proposition}}
\newcommand{\bc}{\begin{corollary}}
\newcommand{\ec}{\end{corollary}}

\newtheorem{lemma}{Lemma}[section]
\newtheorem{corollary}[lemma]{Corollary}
\newtheorem{theorem}[lemma]{Theorem}
\newtheorem{proposition}[lemma]{Proposition}
\newtheorem{definition}[lemma]{Definition}

\newtheorem{example}[lemma]{Example}

\advance\headheight1.15pt
\begin{document}
\title{Fiber cones of ideals with almost minimal multiplicity}

\author[A. V. Jayanthan]{A. V. Jayanthan$^*$}
\author{J. K. Verma}
\address{Department of Mathematics, Indian Institute of Technology
Bombay, Powai, Mumbai, India - 400076}
\email{jayan@math.iitb.ac.in}
\email{jkv@math.iitb.ac.in}
\thanks{$^*$Supported by the National Board for
Higher Mathematics, India}
\thanks{AMS Classification 2000: 13H10, 13H15 (Primary), 13C15, 13A02
(Secondary)}
\maketitle
\begin{abstract}
Fiber cones of 0-dimensional ideals with almost minimal multiplicity 
in Cohen-Macaulay local rings are studied. Ratliff-Rush closure of
filtration of ideals with respect to another ideal is introduced. This
is used to find a bound on the reduction number with respect to an
ideal. Rossi's bound on reduction number in terms of Hilbert
coefficients is obtained as a consequence.
Sufficient conditions are provided for the fiber cone of 0-dimensional
ideals to have almost maximal depth. Hilbert series of such fiber
cones are also computed.
\end{abstract}
\section{Introduction}
Let $(R,\m)$ be a Cohen-Macaulay local ring having infinite residue
field. Let $I$ be an $\m$-primary ideal of $R$ and $K$ an ideal
containing $I$. The fiber cone of $I$ with respect to $K$ is the
standard graded algebra $F_K(I) = \oplus_{n\geq0}I^n/KI^n$. The graded
algebra $F_K(I)$ for $K = \m$ is called the fiber cone $F(I)$ of
$I$. For $K = I$, $F_K(I) = G(I)$, the associated graded ring of $I$.
The objective of this paper is to study the depth of the ring $F_K(I)$
subject to certain conditions imposed on the coefficients of the
Hilbert polynomial $P(I,n)$ corresponding to the Hilbert function
$H(I,n) = \lambda(R/I^n)$, where $\lambda$ denotes the length
function. This theme has been studied for the associated graded rings
by Elias \cite{e}, Goto \cite{g}, Huckaba \cite{h}, Huckaba-Marley
\cite{hm}, Jayanthan-Singh-Verma \cite{jsv}, Rossi \cite{r},
Rossi-Valla \cite{rv}, \cite{rv1}, Sally \cite{s5}, \cite{s2}, Wang
\cite{w}, \cite{w1} and for the fiber cones by Cortadellas-Zarzuela
\cite{cz}, D'Cruz-Raghavan-Verma \cite{crv}, D'Cruz-Verma \cite{dv},
Jayanthan-Verma \cite{jv}, Shah \cite{sh} et al. 

The form ring $G(I)$ and the fiber cone $F(I)$ have been studied
separately. By studying the ring $F_K(I)$ we hope to unify the results
obtained for $G(I)$ and $F(I)$. An initial motivation for this paper
was to find conditions on Hilbert coefficients which will ensure high
depth for $F_K(I)$, thereby obtaining the results known for $G(I)$ and
$F(I)$. We have not been completely successful in providing a unified
approach. However, techniques developed to study the two rings
separately can be unified in the hope of obtaining results for
$F_K(I)$. 

J. D.  Sally studied the Cohen-Macaulay local rings with almost
minimal multiplicity (almost maximal embedding dimension), i.e.,
$\mu(\m) = e_0(\m) + d - 2$. She conjectured that the associated
graded rings of the maximal ideal of such rings have almost maximal
depth, i.e., grade, $\gamma(I)$, of the maximal homogeneous ideal of
$G(I)$ is at least $d-1$. This conjecture was settled in the
affirmative independently by M.  E. Rossi and G. Valla \cite{rv} and
H. -J. Wang \cite{w}. Later M. E.  Rossi extended this result to
$\m$-primary ideals in \cite{r}. She proved that if $I$ is an
$\m$-primary ideal with a minimal reduction $J$ in a Cohen-Macaulay
local ring such that $\lambda(I^2/JI) = 1$, then the associated graded
ring $G(I)$ has almost maximal depth. In section 4 we prove the main
result of the paper, Theorem \ref{ammfiber}.  We prove that if $I$ is
an $\m$-primary ideal with almost minimal multiplicity (i.e.,
$\lambda(\m I/\m J) = 1$ for any minimal reduction $J$ of $I$) and
$\gamma(I) \geq d-2$, then $F(I)$ has almost maximal depth. The method
of the proof is inspired by the methods employed by M. E. Rossi in
\cite{r}. The first main ingredient of the proof is a bound on the
$K$-reduction number, $r_J^K(I) = \min \{n \mid KI^{n+1} = KJI^n \}$.
Such a bound for the usual reduction number $r(I)$ was provided by M.
E. Rossi \cite{r}. This bound played a crucial role in the solution
to an analogue of Sally's conjecture for $\m$-primary ideals by M. E.
Rossi. By specializing the bound on the reduction number $r_J^K(I)$
for $K = \m$, we are able to use it for the fiber cone $F(I)$. 

The second main ingredient of the main theorem is the notion of
Ratliff-Rush closure, $rr_K(I_n) = \cup_{n\geq0}KI_{n+k} : I_1^k$, of
a filtration, $\F = \{I_n\}_{n\geq0}$ of ideals, with respect to an
ideal $K$ containing $I_1$. We will develop the basic properties of
$rr_K(I_n)$ in section 2. We shall find an analogue of Huneke's
fundamental lemma \cite{h2} for the Hilbert function
$\lambda(R/rr_K(I_n))$. As a consequence of this generalization, we
shall provide formulas, in dimension 2, for the coefficients of the
Hilbert polynomial, $P_K(\F,n)$ corresponding to the Hilbert function
$H_K(\F,n) = \lambda(R/KI_n)$. These formulas are crucial for
obtaining the bound on the reduction number $r_J^K(I)$ in Corollary
\ref{redbound}. We shall recover Rossi's bound \cite{r} for $r(I)$ as
a consequence of our bound for $r_J^K(I)$.

One of the motivations for finding numerical conditions which ensure
high depth for $G(I)$ and $F(I)$ is to compute the Hilbert series.
Because of high depth one can work in dimension 1 or 2 where
computation of Hilbert series is relatively easy. As a result, by
imposing conditions on the multiplicity and minimum number of
generators, one can predict the Hilbert series.  In the final section
of this paper, we obtain a formula for the generating function,
$\sum_{n\geq0}H_K(I,n)t^n$, where $I$ is an $\m$-primary ideal with
almost minimal multiplicity with respect to $K$. This formula
generalizes results of Sally and Rossi-Valla. 

\section{Ratliff-Rush closure of a filtration of ideals with respect
to an ideal}

Let $(R,\m)$ a Noetherian local ring of dimension $d > 0$. A
filtration of ideals $\F = \{I_n\}_{n \geq 0}$ is said to be an
$I_1$-good filtration if 
\begin{enumerate}
\item[(i)] $I_nI_m \subseteq I_{n+m}$ for all $n, m \geq 0$,
\item[(ii)] there exists $p \geq 0$ such that $I_1I_n = I_{n+1}$ for
all $n \geq p$ and
\item[(iii)] $I_1$ is $\m$-primary.
\end{enumerate}

For a Hilbert filtration $\F$, let
$H(\F,n) = \lambda(R/I_n)$ denote the Hilbert function of $\F$ and
$$
P(\F,n) = e_0(\F){n+d-1 \choose d} - e_1(\F){n+d-2 \choose d-1} +
\cdots + (-1)^de_d(\F)
$$
denote the corresponding polynomial. 

Let $K$ be an ideal such that $I_{n+1} \subseteq KI_n$ for all $n \geq
0$. Let $H_K(\F,n) = \lambda(R/KI_n)$ be the Hilbert function of $\F$
with respect to $K$. Since $H_K(\F,n) = \lambda(R/I_n) +
\lambda(I_n/KI_n)$, $H_K(\F,n)$ coincides with a polynomial for $n \gg
0$. Let the corresponding polynomial be denoted by
$$
P_K(\F,n) = g_0{n+d-1 \choose d} - g_1{n+d-2 \choose d-1} +
\cdots + (-1)^dg_d
$$

Let $(R,\m)$ be a Noetherian local ring and $\F$ an $I_1$-good
filtration of $R$. Let $K$ be an ideal of $R$ such that $I_1 \subseteq
K$. Let $F_K(\F) = \oplus_{n\geq0} I_n/KI_n$ be the fiber cone of $\F$
with respect to $K$. For $x \in I_1 \backslash KI_1$, let $x^*$ denote
its initial form in the associated graded ring, $G(\F) =
\oplus_{n\geq0}I_n/I_{n+1}$, of $\F$ and $x^o$ denote its initial form
in the fiber cone $F_K(\F)$.
\vskip 3mm
\noindent

We begin by recalling some of the properties of superficial elements
in $F_K(\F)$ proved in \cite{jv}.

\begin{proposition}\label{supeltsprprties}
Let $(R,\m)$ be a Noetherian local ring of dimension $d$ with $R/\m$
infinite. Let $\F = \{I_n\}$ be a Hilbert filtration of $R$, $K$ an
ideal such that $I_{n+1} \subseteq KI_n$ for all $n \geq 0$. Then
\begin{enumerate}
\item[1.] There exists an $x \in I_1 \backslash KI_1$ such that $x^o$ is 
superficial in $F_K(\F)$ and $x^*$ is superficial in $G(\F)$.

\item[2.] If, for $x \in I_1 \backslash KI_1$, $x^o$ is superficial in
$F_K(\F)$ and $x^*$ is superficial in $G(\F)$, then there exists a $c 
> 0$ such that $(KI_n : x) \cap I_c = KI_{n-1}$ for all $n > c$.
Moreover if $x$ is regular in $R$, then $KI_n : x = KI_{n-1}$ for all
$n \gg 0$. 

\item[3.] If $x \in I_1 \backslash \; KI_1$ is regular in R, $x^o$ is
regular in $F_K(\F)$ and $x^*$ is regular in $G(\F)$, then $KI_n : x =
KI_{n-1}$ for all $n \geq 1$.

\item[4.] Let $x \in I_1$ be such that $x^*$ is superficial in $G(\F)$ and
$x^o \in F_K(\F)$ is superficial in $F_K(\F)$. Let $\bar{\F} = \{I_n +
xR/xR\}_{n\geq 0}$ and $\bar{K} = K/xR$. If $\depth
F_{\bar{K}}(\bar{\F}) > 0$, then $x^o$ is regular in $F_K(\F)$.

\item[5.] Let $x_1, \ldots, x_k \in I_1$. Assume that 
\begin{enumerate}
\item[$(i)$]  $x_1, \ldots, x_k$ is a regular sequence in $R$.
\item[$(ii)$] $x_1^* \ldots, x_{k}^* \in G(\F)$ is a regular sequence.
\item[$(iii)$] $x_1^o, \ldots, x_k^o \in F_K(\F)$ is a superficial 
sequence.
\end{enumerate}
Then $\depth_{(x_1^o, \ldots, x_k^o)}F_K(\F) = k$ if and only if
$(x_1, \ldots, x_k) \cap KI_n = (x_1, \ldots, x_k)KI_{n-1}$ for
all $n \geq 1$.
\end{enumerate}
\end{proposition}

\begin{definition} The Ratliff-Rush closure of $\F = \{I_n\}$ with
respect to $K$ is the sequence of ideals $rr_K(\F) = \{rr_K(I_n)\}_{n
\geq 0}$ given by 
$$
rr_K(I_n) = \bigcup_{k \geq 1} (KI_{n+k} : I_1^k). 
$$
\end{definition}
\noindent
The Ratliff-Rush closure of a filtration of ideals with respect to an
ideal behaves quite similar to the Ratliff-Rush closure of an ideal.
We summarize some of its properties.

\begin{proposition}\label{rrproperties}
\begin{enumerate}
\item[1.] $rr_K(I_n) = \bigcup_{k\geq1}(KI_{nk+n} : I_n^k)$.
\item[2.] If $\grade I_1 > 0$, then $rr_K(I_n) = KI_n$ for $n \gg 0$.
\item[3.] If $J = (x_1, \ldots, x_s)$ is a reduction of $I_1$, then
$rr_K(I_n) = \bigcup_{k \geq 1}(KI_{n+k} : (x_1^k, \ldots, x_s^k)).$
\item[4.] If $J$ is a reduction of $I_1$, then
$rr_K(I_n) = \bigcup_{k \geq 1}(KI_{n+k} : J^k)$.
\item[5.] If $J$ is a minimal reduction of $I_1$, then 
$rr_{K}(I_{n}) : J = rr_{K}(I_{n-1})$ for all $n \geq 1$.
\end{enumerate}
\end{proposition}
\begin{proof}
1. 
Note that $KI_{n+1} : I_1 \subseteq KI_{n+2} : I_1^2 \subseteq \cdots$
is an increasing chain of ideals in $R$. Hence $rr_K(I_n) = KI_{n+k} :
I_1^k$ for $k \gg 0$. Since the chain $KI_{n+n} : I_n \subseteq
KI_{n+2n} : I_n^2 \subseteq \cdots$ also terminates, it is
enough to show that $rr_K(I_n) = KI_{nk+n} : I_n^k$ for $k \gg 0$.
Suppose $x \in KI_{nk+n} : I_n^k$. 
Since $I_1^{nk} \subseteq I_n^k$, $xI_1^{nk} \subseteq KI_{nk+n}$.
Therefore for $k \gg 0$, $x \in KI_{nk+n} : I_1^{nk} = rr_K(I_n)$.
Conversely, let $x \in rr_K(I_n)$. Since $\F$
is an $I_1$-good filtration, there exists $p_0$ such that  $I_1I_p = I_{p+1}$
for $p \geq p_0$. Choose $k \gg 0$. Then 
$$
xI_n^k \subseteq xI_{nk} \subseteq xI_{p_0}I_1^{nk-p_0} \subseteq
I_{p_0}KI_{nk-p_0+n} \subseteq KI_{nk+n}.
$$
Therefore $x \in KI_{nk+n} : I_n^k$, so that $rr_K(I_n)  = KI_{nk+n} :
I_n^k$ for $k \gg 0$.
\vskip 2mm
\noindent
2. Let $x \in I_1$ be such that $x$ is regular in $R$ and $x^o$ is
superficial in $F_K(\F)$ and $x^*$ is superficial in $G(\F)$. Then, by
Proposition \ref{supeltsprprties}(2), $KI_n : x = KI_{n-1}$ for $n \gg 
0$. Therefore
$KI_n \subseteq KI_{n+1} : I_1 \subseteq KI_{n+1} : x = KI_n$ for $n
\gg 0.$ Thus $KI_{n+1} : I_1 = KI_n$ for $n \gg 0$. We show that
$KI_{n+k} : I_1^k = KI_n$ for all $k \geq 1$. Apply induction on $k$.
The result is proved for $k = 1$. Assume that the result is true for
$k - 1$. Then
\begin{eqnarray*}
KI_{n+k} : I_1^k & = & (KI_{n+k} : I_1^{k-1}) : I_1 \\
& = & KI_{n+1} : I_1 \hspace*{3mm}\mbox{(by induction)}\\
& = & KI_n \hspace*{2mm} \mbox{ for } n \gg 0.
\end{eqnarray*}
Therefore $rr_K(I_n) = KI_n$ for $n \gg 0$.
\vskip 2mm
\noindent
3. Let $(\underline{x}) = (x_1, \ldots, x_s)$ and $(\underline{x})^{[k]}
= (x_1^k, \ldots, x_s^k).$ Clearly $KI_{n+k} : I_1^k \subseteq 
KI_{n+k} : (\underline{x})^{[k]}$. Since $(\underline{x})$ is a
reduction of $I_1$, there exists an integer $r$ such that
$(\underline{x})^mI_1^n = I_1^{n+m}$ for all $n \geq r$ and $m \geq 1$.
Let $z \in KI_{n+k} : (\underline{x})^{[k]}$ for $k \gg 0$. Then 
\begin{eqnarray*}
zI_1^{r+sk} & = & z(\underline{x})^{sk}I_1^r \\
& = &  \left(\sum_{\mid \alpha \mid = sk}z x_1^{\alpha_1}\cdots
x_s^{\alpha_s}\right)I_1^r \\
& \subseteq & \sum_{\mid \alpha \mid =
sk}KI_{n+\alpha_i}x_1^{\alpha_1} \cdots \widehat{x_i^{\alpha_i}}\cdots
x_s^{\alpha_s}I_1^r \hspace*{3mm}\mbox{ where }\alpha_i \geq k \\
& \subseteq & KI_{n+r+sk}.
\end{eqnarray*}
Therefore $z \in KI_{n+r+sk} : I_1^{sk+r} = rr_K(I_n)$. 
\vskip 2mm
\noindent
4. Let $J = (x_1, \ldots, x_s)$. Since $KI_{n+k} : I_1^k \subseteq
KI_{n+k} : J^k \subseteq KI_{n+k} : (x_1^k, \ldots, x_s^k)$ for all
$k$, the assertion follows.
\vskip 2mm
\noindent
5. For $k \gg 0$, we have 
\begin{eqnarray*}
rr_{K}(I_n) : J & = &(K I_{n+k} : J^{k}) : J \\
& = & K I_{n+k} : J^{k+1} = rr_{K}(I_{n-1}).
\end{eqnarray*}
\end{proof}
\vskip 2mm
\noindent
The next lemma inspired the definition of Ratliff-Rush closure of a
filtration of ideals with respect to another ideal. 
\vskip 2mm
\noindent
\begin{proposition}\label{0locfiber}
Let $F=F_K(\F)$ and let
$[H^0_{F_+}(F)]_n$ denote the $n$-th graded component of
the local cohomology module $H^0_{F_+}(F)$.  Suppose $\grade I_1 > 0$.
Then for all $n \geq 0$,
$$
[H^0_{F_+}(F)]_n = \frac{rr_K(I_n) \cap I_n}{KI_n}.
$$
If $\grade (I_1) > 0$ and $\gamma(\F) > 0$, then 
$\grade F_+ > 0$ if and only if $rr_K(I_n) = KI_n$ for all $n \geq 0$.
\end{proposition}
\begin{proof}
Let $y \in I_n$ and $y^o \in [H^0_{F_+}(F)]_n = 0 :_{F_n} F_+^k$ for $k
\gg 0$. Then $y^o F_+^k = 0$. Therefore $yI_1^k \subseteq KI_{n+k}$.
Hence $y \in (KI_{n+k} : I_1^k) \cap I_n = 
rr_K(I_n) \cap I_n$. Therefore $[H^0_{F_+}(F)]_n \subseteq (rr_K(I_n) 
\cap I_n)/KI_n$. Suppose $y^o \in (I_n \cap rr_K(I_n))/KI_n$. Then,
for $k \gg 0$, $yI_1^k \subseteq KI_{n+k}$. Therefore $y^o F_+^k = 0$
so that $y^0 \in 0 :_{F_n} F_+^k \subseteq [H^0_{F_+}(F)]_n$.
\vskip 2mm
\noindent
Suppose $rr_K(I_n) = KI_n$ for all $n \geq 0$. Then $H^0_{F_+}(F) =
\oplus_{n \geq 0} I_n \cap rr_K(I_n) /KI_n = 0$. Therefore $\grade F_+
> 0$. Conversely, suppose $\grade F_+ > 0$. Then $rr_K(I_n) \cap I_n =
KI_n$ for all $n \geq 0$. Suppose $y \in rr_K(I_n) = KI_{n+k} :
I_1^k$. Choose a regular element $x_1 \in I_1$ such that 
$x_1^o$ is regular in $F$ and $x_1^*$ is regular in $G(\F)$.
Then $yx_1^k \in KI_{n+k}$ so that $yx_1^k \in KI_{n+k} \cap
(x_1^k) = x_1^kKI_n$, by Lemma \ref{supeltsprprties}(3).  Therefore 
$y \in KI_n$ and hence $rr_K(I_n) = KI_n$. 
\end{proof}
\vskip 2mm
\noindent
In the next proposition we obtain a generalization of Huneke's 
fundamental lemma \cite{h2} for the function $H_K(\F,n)$. It also
shows that once we know a minimal reduction of $I_1$, we can compute
the coefficients $g_1$ and $g_2$ and hence the Hilbert polynomial of
$F_K(\F)$ can be completely  determined.
\begin{proposition}\label{genfunction1}
Let $(R,\m)$ be a $2$-dimensional Cohen-Macaulay local ring. Let $\F =
\{I_n\}$ be a Hilbert filtration of $R$. Let $J = (x,y)$ be a minimal
reduction of $I_1$.  Then  $\lambda(R/rr_K(I_n))$ coincides with the 
polynomial $P_{K}(\F,n)$, for $n \gg 0$ and the following are true: \\
$1$. For $n \geq 2$, 
$$
\Delta^2[P_{K}(\F, n) - \lambda(R/rr_K(I_n))] =
\lambda\left(\frac{rr_{K}(I_n)}{Jrr_{K}(I_{n-1})}\right). 
$$
$2$. Set 
$$
v_n = \left\{
  \begin{array}{lll}
   e_0(\F) - \lambda(R/rr_K(I_0)) \mbox{ if } n = 0 \\
   e_0(\F) - \lambda(R/rr_K(I_1)) + 2 \lambda(R/rr_K(I_0)) \mbox{ if }
   n = 1 \\
   \lambda\left(\frac{rr_{K}(I_n)}{Jrr_{K}(I_{n-1})}\right) \mbox{ if }
   n \geq 2. 
  \end{array}
 \right.
$$
Then, $g_1 = \sum_{n \geq 1}v_n - \lambda(R/rr_K(I_0))$ and
$g_2 = \sum_{n\geq1}(n-1)v_n - \lambda(R/rr_K(I_0)).$ 

\end{proposition}
\begin{proof}
1. Since $rr_{K}(I_n) = K I_n$ for $n \gg 0$,
$\lambda(R/rr_K(I_n)) = H_{K}(\F,n) = P_{K}(\F,n)$ for $n \gg 0$.
Consider the exact sequence:
$$
\begin{CD}
0 @>>> \frac{R}{rr_K(I_{n-1}) : J} @>\beta>> 
\left(\frac{R}{rr_K(I_{n-1})}\right)^2
@>\alpha>> \frac{J}{Jrr_K(I_{n-1})} @>>> 0,
\end{CD}
$$
where the maps $\alpha$ and $\beta$ are defined as, $\alpha(\bar{r}, 
\bar{s}) = \overline{xr + ys}$ and $\beta(\bar{r}) = (\bar{y}\bar{r},
-\bar{x}\bar{r})$.  It follows that for all $n \geq 2$,
\begin{eqnarray*}
2\lambda(R/rr_K(I_{n-1})) & = & \lambda(R/(rr_K(I_{n-1}) : J)) +
\lambda(J/Jrr_K(I_{n-1})) \\
& = & \lambda(R/(rr_K(I_{n-1}) : J)) + \lambda(R/Jrr_K(I_{n-1})) -
\lambda(R/J).
\end{eqnarray*}
Therefore $e_0(\F) + 2\lambda(R/rr_K(I_{n-1})) = \lambda(R/Jrr_K(I_{n-1})) +
\lambda(R/(rr_K(I_{n-1}) : J))$. Hence
\begin{eqnarray*}
e_0(\F) & - & \lambda(R/rr_K(I_{n})) + 2\lambda(R/rr_K(I_{n-1})) - 
\lambda(R/rr_K(I_{n-2}))\\
& = & \lambda(R/Jrr_K(I_{n-1})) - \lambda(R/rr_K(I_{n})) + 
\lambda(R/(rr_K(I_{n-1}) :
J))- \lambda(R/rr_K(I_{n-2})) \\
& = & \lambda(rr_K(I_n)/Jrr_K(I_{n-1})) - 
\lambda(rr_K(I_{n-1}) : J/rr_K(I_{n-2}))
\end{eqnarray*}
Since $\Delta^2P_K(\F,n) = e_0(\F)$, 
$$
\Delta^2\left[P_K(\F,n) - H_K(\F,n)\right] =
\lambda\left(\frac{rr_K(I_n)}{Jrr_K(I_{n-1})}\right) - 
\lambda\left(\frac{rr_K(I_{n-1}) : J}{rr_K(I_{n-2})}\right).
$$
By Proposition \ref{rrproperties}(5), we have $rr_K(I_n) : J =
rr_K(I_{n-1})$ for all $n \geq 1$. 
Therefore for all $n \geq 2,$
$$
\Delta^2[P_{K}(\F, n) - \lambda(R/rr_K(I_n))] =
\lambda\left(\frac{rr_{K}(I_n)}{Jrr_{K}(I_{n-1})}\right).
$$
\vskip 2mm
\noindent
2. Write $P_{K}(\F,n) = e_0(\F){n+2 \choose 2} - g_1'(n+1) + g_2'$.
Then $\sum_{n\geq0}\Delta^2P_K(\F,n)t^n = e_0(\F)/(1-t)$. Let
$\sum_{n\geq0}\lambda(R/rr_K(I_n))t^n = f(t)/(1-t)^3$. Then $g_1' = f'(1)$ 
and $g_2' = f''(2)/2$, by Proposition 4.1.9 of \cite{bh}. Note that
for the Ratliff-Rush filtration, for all $n < 0$,
$\lambda(R/rr_K(I_n)) = 0$. We have,

\begin{eqnarray*}
\sum_{n\geq0}\Delta^2[P_K(\F,n) - \lambda(R/rr_K(I_n))]t^n 
& = & \frac{e_0(\F)}{(1-t)} -
\frac{f(t)}{(1-t)} - 2\lambda(R/rr_K(I_{-1}))\\
& + & \lambda(R/rr_K(I_{-2}))) + t\lambda(R/rr_K(I_{-1}))  \\
 & = & \frac{e_0(\F) - f(t)}{(1-t)}\\
 & = & \sum_{n\geq0}v_nt^n
\end{eqnarray*}
Therefore $e_0(\F) - f(t) = (1-t)\sum_{n\geq0}v_nt^n$ and hence $f(t)
= e_0(\F) - (1-t)\sum_{n\geq0}v_nt^n$. Therefore 
$$
f'(t) = \sum_{n\geq0}v_nt^n - (1-t)\sum_{n\geq0}nv_nt^{n-1} ~~~
\mbox{ and }
$$ 
$$f''(t) =
\sum_{n\geq0}nv_nt^{n-1} - (1-t)\sum_{n\geq0}n(n-1)v_nt^{n-2} +
\sum_{n\geq0}nv_nt^{n-1}.
$$ 
Hence $g_1' = f'(1) = \sum_{n\geq0}v_n$ and
$g_2' = f''(1)/2 = \sum_{n\geq0}nv_n$. Therefore 
$$g_1 = g_1' - e_0(\F) = \sum_{n \geq 1}v_n - \lambda(R/rr_K(I_{0}))~~~
\mbox{ and }
$$
$$
g_2 = g_2' - g_1 = \sum_{n\geq0}nv_n - \sum_{n \geq 1}v_n - 
\lambda(R/rr_K(I_0)) = \sum_{n\geq1}(n-1)v_n - \lambda(R/rr_K(I_0)).
$$
\end{proof}

\vskip 2mm
\noindent
As a consequence we derive formulas for the Hilbert coefficients 
obtained by Huneke for the $I$-adic filtration in \cite{h2}.

\begin{corollary}
Let $(R,\m)$ be a $2$-dimensional Cohen-Macaulay local ring and $\F$
be a Hilbert filtration. Then 
$e_1(\F) = \sum_{n\geq1}v_n$ and $e_2(\F) = \sum_{n\geq1}(n-1)v_n$.
\end{corollary}
\begin{proof}
Put $K = R$ in Proposition \ref{genfunction1}(2) and note that, for $K
= R$, $rr_K(I_0) = R$. 
\end{proof}

\section{Bounds on reduction numbers}

In this section we obtain a bound on the $K$-reduction number of an
$\m$-primary ideal (see Definition \ref{rednK}) from which we derive
Rossi's bound (Corollary 1.5, \cite{r}) for the reduction number. We
use this bound to prove the almost maximal depth condition for the
fiber cone. We set the following notation for the rest of the section.
Let $(R, \m)$ denote a Cohen-Macaulay local ring with infinite residue
field.  Let $I$ be an $\m$-primary ideal of $R$ and let $J$ be its
minimal reduction. Let $K$ be an ideal containing $I$ and let
$rr_{K}(I^n)$ denote the Ratliff-Rush closure of $I^n$ with respect to
$K$. For $n \geq 0$, set 
$$
\rho_n^{K}= \lambda(rr_{K}(I^{n+1})/Jrr_{K}(I^{n})) \mbox{ and }
\nu_n^{K} = \lambda(KI^{n+1}/KJI^{n}).
$$

\vskip 2mm
\noindent
The next proposition, due to M. E. Rossi, played a crucial
role in solving the conjecture of Sally \cite{rv} and its 
generalization to case of $\m$-primary ideals \cite{r}.
For an ideal $I$ in $R$, let $\R(I) = \oplus_{n\geq0}I^nt^n$ denote
the Rees algebra of $I$. For an $\R(I)$-module $M$, 
put $\ann_{I^{\nu}}(M) = \{x \in I^{\nu} \mid xt^{\nu}M = 0\}$.

\begin{proposition}\label{rossi1} Let $I$ be an ideal of a Noetherian
local ring $R$ and let $J$ be a minimal reduction of $I$. Let $M$ be
an $\R(I)$-module of finite length as $R$-module. Let $\nu$ be the
minimum number of generators of $M/\R(J)_+M$ as an $R$-module. Then 
$$
I^{\nu} = JI^{\nu-1} + \ann_{I^{\nu}}(M).
$$
\end{proposition}

\begin{definition}\label{rednK}
Let $J$ be a minimal reduction of an ideal $I$. Put
$$
r_J^{K}(I) := \min \{n \mid KI^{n+1} = KJI^n\}.
$$ 
The integer $r_J^{K}(I)$ is called the $K$-reduction number of $I$ 
with respect to $J$.
\end{definition}
We now give a bound for the $K$-reduction number of an $\m$-primary 
ideal. 

\begin{theorem}\label{rednobound}
Let $(R,\m)$ be a Cohen-Macaulay local ring of dimension $d > 0$. 
Let $I$ be an $\m$-primary ideal of $R$ and let $J$ be
a minimal reduction of $I$. Let $K$ be an ideal containing $I$ such
that $KI \cap J = KJ$. Then 
$$
r_J^{K}(I) \leq
\sum_{j\geq0}\rho_j^K - \lambda(KI/KJ) + 1.
$$
\end{theorem}

\begin{proof}
Let $M := \bigoplus_{n\geq1}rr_{K}(I^n)/KI^n$. Then $M$ is a
finitely generated $\R(I)$-module and $\lambda_R(M) < \infty$, by
\ref{rrproperties}(2).  For $j \geq 0$, $[M/\R(J)_+M]_{j+1} =
M_{j+1}/(J^{j+1}M_0+J^jM_1+ \cdots + JM_j)$. For $1 \leq i \leq j+1$
and $k \gg 0$, we have
\begin{eqnarray*}
J^iM_{j-i+1} & = & J J^{i-1}M_{j-i+1} \subseteq JI^{i-1}M_{j-i+1} \\
&=& \frac{JI^{i-1}(KI^{j+1-i+k} : I^k) + KI^{j+1}}{KI^{j+1}} \\
&\subseteq & \frac{J rr_{K}(I^j) + KI^{j+1}}{KI^{j+1}} \\
& = & JM_j.
\end{eqnarray*}
Therefore
$[M/\R(J)_+M]_{j+1} \cong rr_{K}(I^{j+1})/J
rr_{K}(I^j) + KI^{j+1}.$ We have
$$
\lambda\left(rr_{K}(I^{j+1})/Jrr_{K}(I^j)+KI^{j+1}\right) \leq 
\lambda\left(rr_{K}(I^{j+1})/Jrr_{K}(I^j)\right)
$$
and equality occurs if and only if $KI^{j+1} \subseteq Jrr_{K}(I^j)$.
Since $J$ is a reduction of $I$ and $rr_{K}(I^n) = KI^n$ for $n \gg
0$, there exists a $j$ such that $KI^{j+1} = KJI^j = J rr_{K}(I^j)$.
Let $k = \min\{j \mid KI^{j+1} \subseteq J rr_{K}(I^j) \}.$ Let
$\mu_j$ be the minimum number of generators of
$rr_{K}(I^{j+1})/Jrr_{K}(I^j)+KI^{j+1}$ as an $R$-module. Then $\mu_j
\leq \lambda(rr_{K}(I^{j+1})/Jrr_{K}(I^j)+ KI^{j+1})$. Let $\mu =
\sum_{j\geq0}\mu_j$. Then by Proposition \ref{rossi1}, $I^{\mu} =
JI^{\mu-1} + \ann_{I^{\mu}}(M)$. Therefore
\begin{eqnarray*}
KI^{\mu+k+1} & = & KI^{k+1}(JI^{\mu-1}+\ann_{I^{\mu}}(M)) \\
& = & KJI^{\mu+k} + KI^{k+1}\ann_{I^{\mu}}(M) \\
& \subseteq & KJI^{\mu+k} + Jrr_{K}(I^k)\ann_{I^{\mu}}(M).
\end{eqnarray*}
Since $J rr_{K}(I^k)\ann_{I^{\mu}}(M) \subseteq KJI^{\mu+k}$, 
$KI^{\mu+k+1} = KJI^{\mu+k}$. Therefore 
\begin{eqnarray*}
r_J^{K}(I) & \leq & \mu + k = \sum_{j\geq0} \mu_j + k \leq
\sum_{j\geq0} \lambda\left(\frac{rr_{K}(I^{j+1})}{Jrr_{K}(I^j) +
KI^{j+1}} \right) + k \\
& = & \lambda\left(rr_K(I)/Jrr_K(I^0)+KI\right) + \sum_{j\geq1}\left[\rho^K_j -
\lambda(Jrr_K(I^j)+KI^{j+1}/Jrr_K(I^j))\right] + k.
\end{eqnarray*}

Since $KI \cap J = KJ$, $KJ \subseteq Jrr_K(I^0)\cap KI \subseteq
J\cap KI =KJ$ and hence 
$$
\lambda\left(rr_K(I)/Jrr_K(I^0)+KI\right) = \rho^K_0 - \lambda(KI/KJ).
$$ 
Also note that, 
$$
\lambda\left(\frac{rr_{K}(I^{j+1})}{Jrr_{K}(I^j) + KI^{j+1}} \right)
   \leq \left\{
\begin{array}{ll}
\rho_j^K - 1 \mbox{ for
} j = 1, \ldots, k-1 ~~\mbox{ and }\\
\rho_j^K \mbox{ for } j \geq k.
\end{array}\right.
$$

Therefore we get,
\begin{eqnarray*}
r^K_J(I) & \leq & \rho_0^K - \lambda(KI/KJ) + \sum_{j=1}^{k-1} [\rho_j^K - 1] 
+ \sum_{j\geq k} \rho_j^K + k \\
 & = & \sum_{j\geq0}\rho_j^K - \lambda(KI/KJ) + 1.
\end{eqnarray*}
\end{proof}

\vskip 2mm
\noindent
The following lemma is quite well-known. We include it for the sake of 
completeness.
\begin{lemma}\label{genlemma1}
Let $(R,\m)$ be a Noetherian local ring and let $J=(x_1, \ldots, x_s)$
be an ideal generated by a regular sequence in $R$. Then for any ideal
$K$ containing $J$, $J/KJ \cong (R/K)^s$. 
\end{lemma}
\begin{proof}
Consider the map $\phi : (R/K)^s \longrightarrow J/KJ$, defined as 
$$
\phi(\bar{r}_1, \ldots, \bar{r}_s) = \overline{r_1x_1+\cdots+r_sx_s}.
$$
The map $\phi$ is clearly surjective. Suppose for some $r_1, \ldots,
r_s \in R$, $r_1x_1+\cdots+r_sx_s \in KJ$. Write $r_1x_1+\cdots+r_sx_s
= t_1x_1+\cdots+t_sx_s$ for some $t_1, \ldots, t_s \in K$. Then
$(r_1-t_1)x_1 = (t_2 -r_2)x_2+\cdots+(t_s-r_s)x_s$. Since $x_1,
\ldots, x_s$ is a regular sequence, $r_1-t_1 \in (x_2, \ldots, x_s)
\subseteq K$ and hence $r_1 \in K$. Similarly $r_i \in K$ for all $i =
1, \ldots, s$. Therefore $\phi$ is an isomorphism.
\end{proof}

We obtain a bound on the reduction number $r^K_J(I)$ in terms of the 
Hilbert coefficient $g_1$.
\vskip 2mm
\noindent

\begin{corollary}\label{redbound}
Let $(R,\m)$ be a $2$-dimensional Cohen-Macaulay local ring, $I$ an
$\m$-primary ideal, $K$ an ideal containing $I$ and $J$ a minimal
reduction of $I$. If $KI \cap J = KJ$, then
$$
r_J^{K}(I) \leq g_1 - \lambda(KI/KJ) + 1 + \lambda(R/K).
$$
\end{corollary}
\begin{proof}
By Theorem \ref{rednobound} we get 
\begin{eqnarray*}
r_J^{K}(I) & \leq & \lambda(rr_K(I)/Jrr_K(I^0)) + 
\sum_{j\geq1}\rho_j^K + 1 - \lambda(KI/KJ) \\
& = & \lambda(R/Jrr_K(I^0)) -
\lambda(R/rr_K(I)) + \sum_{j\geq1}\rho_j^K + 1 -
\lambda(KI/KJ) \\ 
& = & \lambda(R/J) + \lambda(J/Jrr_K(I^0)) - \lambda(R/rr_K(I)) +
\sum_{j\geq1}\rho_j^K + 1 - \lambda(KI/KJ) \\
& = & e_0(I) + 2\lambda(R/rr_K(I^0)) - \lambda(R/rr_K(I)) + 
\sum_{j\geq1}\rho_j^K + 1 - \lambda(KI/KJ)\\
& \leq &
e_0(I) + \lambda(R/rr_K(I^0)) + \lambda(R/K) - \lambda(R/rr_K(I))
+ \sum_{j\geq1}\rho_j^K + 1 - \lambda(KI/KJ). 
\end{eqnarray*}
The last equality follows from Lemma \ref{genlemma1} and the 
inequality follows since $K \subseteq rr_K(I^0)$. By Lemma
\ref{genfunction1}, $g_1 = e_0(I) + \lambda(R/rr_K(I^0)) -
\lambda(R/rr_K(I)) + \sum_{j\geq1}\rho_j^K$. Therefore
$$
r_J^{K}(I) \leq g_1 - \lambda(KI/KJ) + 1 + \lambda(R/K).
$$
\end{proof}

\begin{corollary}(Rossi's bound)
Let $(R,\m)$ be a Cohen-Macaulay local ring of dimension $2$. Let $I$
be an $\m$-primary ideal of $R$ and $J$ be a minimal reduction of $I$.
Then 
$$
r_J(I) \leq e_1(I) - e_0(I) + \lambda(R/I) + 1.
$$
\end{corollary}
\begin{proof} Put $K = R$ in Corollary \ref{redbound} and note that 
for $K = R$, $g_i = e_i$ for $i = 0, \ldots, d$.
\end{proof}
\vskip 2mm
\noindent
Our objective in introducing $r_J^K(I)$ is to obtain bounds for
$r^{\m}_J(I)$ which in turn is used to study the depth of fiber cones
of ideals with almost minimal multiplicity.
\begin{corollary}\label{ammredbound}
Let $(R,\m)$ be a $2$-dimensional Cohen-Macaulay local ring, $I$ an
$\m$-primary ideal and $J$ a minimal reduction of $I$. Then
$$
r^{\m}_J(I) \leq g_1 + 2 - \lambda(\m I/\m J).
$$
\end{corollary}
\begin{proof}
Put $K = \m$ and use the fact that $\m I \cap J  = \m J$.
\end{proof}
\vskip 2mm
\noindent
If $I$ is an $\m$-primary ideal, $J$ is a minimal
reduction of $I$ and $x^*$ is  regular in $G(I)$, then $r_J(I) =
r_{\bar{J}}(\bar{I})$, \cite{h1}. In the following lemma we prove that
a similar result holds for the $K$-reduction number also.
\begin{lemma}\label{reductionmodulo}
Let $(R,\m)$ be a Noetherian local ring of dimension $d > 0$. Let $I$ 
be an $\m$-primary ideal of $R$, $K$ an ideal containing $I$ and $J$ a
minimal reduction of $I$. Let $x \in I \backslash KI$ be such that
$x^*$ is regular in $G(I)$ and $x^o$ is regular in $F_K(I)$. Then 
$r_{\bar{J}}^{\bar{K}}(\bar{I}) = r_J^K(I)$, where $``-"$ denote
images $\modulo (x)$. 
\end{lemma}
\begin{proof} Clearly $r_{\bar{J}}^{\bar{K}}(\bar{I}) \leq r_J^K(I)$.
Suppose for some $n$, $\bar{K}\bar{I}^n =
\bar{K}\bar{J}\bar{I}^{n-1}$. Then $KI^n + xR = KJI^{n-1} + xR$ and
hence $KI^n = KI^n \cap (KJI^{n-1} + xR) = KJI^{n-1} + (xR \cap
KI^n)$. Since $x^*$ is regular in $G(I)$ and $x^o$ is regular in
$F_K(I)$, by Proposition \ref{supeltsprprties}(5), $xR \cap KI^n=
xKI^{n-1}$.  Hence $KI^n = KJI^{n-1}$. Therefore
$r_{\bar{J}}^{\bar{K}}(\bar{I}) = r_J^{K}(I)$.
\end{proof}

\section{Ideals with almost minimal multiplicity}

Let $(R,\m)$ be a Cohen-Macaulay local ring of dimension $d > 0$. 
Let $I$ be an $\m$-primary ideal and $J$ a minimal reduction of $I$.
Then $\mu(I) + \lambda(\m I/\m J) = e_0(I) - \lambda(R/I) + d$, where
$\mu(I)$ is the minimum number of generators of $I$. Hence $\mu(I)
\leq e_0(I) - \lambda(R/I) + d$ and the equality occurs if and only if
$\m I = \m J$. S. Goto defined an ideal to have minimal multiplicity
if $\m I = \m J$, \cite{g}. He characterized various properties of
the associated graded rings, the Rees algebras and fiber cones of such
ideals. Inspired by these results, we define ideals of almost minimal
multiplicity.

\vskip 2mm
\noindent
\begin{definition}\label{ammdef}
An ideal $I$ is said to have almost minimal multiplicity with respect
to an ideal  $K \supseteq I$ if for any minimal reduction $J$ of
$I$, $\lambda(KI/KJ) = 1$. We say that $I$ has almost minimal
multiplicity if $\lambda(\m I/\m J) = 1$.
\end{definition}
\vskip 1mm
\noindent
{\bf Remark :} For any $\m$-primary ideal $I$, an ideal $K \supseteq
I$ and a minimal reduction $J$ of $I$, $\lambda(KI/KJ) = 1$ if and
only if $\lambda(I/KI) = e_0(I) - \lambda(R/I) + \lambda(J/KJ) - 1 =
e_0(I) - \lambda(R/I) + d \lambda(R/K) - 1$, by Lemma \ref{genlemma1}.
Hence the definition of almost minimal multiplicity with respect to 
$K$ is independent of the minimal reduction $J$ chosen for $I$.

For $K = I$, the almost minimal multiplicity condition is equivalent
to $\lambda(I/I^2) = e_0(I) - (d-1) \lambda(R/I) - 1$, which was the 
condition imposed on the ideal in \cite{r} to obtain the almost
maximal depth for the associated graded ring.
\vskip 2mm
\noindent
\begin{lemma}\label{genlemma2}
Let $(R,\m)$ be a Cohen-Macaulay local ring and let $I$ be an ideal
with almost minimal multiplicity with respect to $K \supseteq I$. Then
for any minimal reduction $J$ of $I$, $\lambda(KI^n/KJI^{n-1}) = 1$
for all $n = 1,\ldots, r_J^K(I).$
\end{lemma}
\begin{proof}
Since $I$ has almost minimal multiplicity, $\lambda(KI/KJ) = 1$ for
any minimal reduction $J$ of $I$. Let $a \in K, \; b \in I$ be such
that $KI = KJ + (ab)$ and $\m ab \subseteq KJ$. Then it can easily be
seen by induction that $KI^n = KJI^{n-1} + (ab^n)$ with $\m ab^n
\subseteq KJI^{n-1}$. Hence $\lambda(KI^n/KJI^{n-1}) = 1$ for all
$n = 1,\ldots, r_J^K(I).$
\end{proof}
In this section we show that if $\gamma(I) \geq d-2$, then the fiber
cone of an $\m$-primary ideal $I$ with almost minimal multiplicity has
almost maximal depth. The method of the proof is analogous to the
method employed by M. E. Rossi in \cite{r} to prove the almost maximal
depth condition for the associated graded ring. 

\begin{lemma}\label{redinvariance}
$1.$ Let $(R,\m)$ be a $2$-dimensional Cohen-Macaulay local ring. Let
$I$ be an $\m$-primary ideal of $R$ with almost minimal multiplicity
and let $J$ be a minimal reduction of $I$. Then $r_J^{\m}(I) \leq g_1
+ 1$. 
\vskip 2mm
\noindent
$2.$ Let $x \in I \backslash \m I$ be such that $x^o$ is superficial
in $F(I)$ and $x^*$ is superficial in $G(I)$.  Let $``-"$ denote
images $\modulo (x)$. If $\bar{I}$ has almost minimal multiplicity,
then $r_{\bar{J}}^{\bar{\m}}(\bar{I}) = r_J^{\m}(I) = g_1+1$.
\end{lemma}
\begin{proof}1. The inequality $r_J^{\m}(I) \leq g_1+1$ directly follows
from Corollary \ref{ammredbound}.  
\vskip 2mm
\noindent
2. Set $s = r_{\bar{J}}^{\bar{\m
}}(\bar{I})$. Clearly $s \leq r_J^{\m}(I)$. As $x^o$ is superficial in
$F(I)$ and $x^*$ is superficial in $G(I)$, $g_1 = \bar{g}_1$, where
$\bar{g}_i$ denote coefficients of the polynomial corresponding to
$\lambda(\bar{R}/\bar{\m}\bar{I}^n)$.  Since $\dim \bar{R} = 1$, by
Theorem 5.3 of \cite{jv}, $\bar{g}_1 = \sum_{n\geq1}
\lambda(\bar{\m}\bar{I}^n/\bar{\m}\bar{J}\bar{I}^{n-1}) - 
\lambda(\bar{R}/\bar{\m})$. Since $\bar{I}$ has almost minimal
multiplicity, $\lambda(\bar{\m}\bar{I}^n/\bar{\m}\bar{J}\bar{I}^{n-1})
= 1$ for $n = 1, \ldots, s$. Therefore $g_1 = s- 1$ and hence 
$r_J^{\m}(I) = r_{\bar{J}}^{\bar{\m}}(\bar{I}) = g_1+1$.
\end{proof}

\vskip 2mm
\noindent
We now prove the main result of this section. 

\vskip 2mm
\noindent
\begin{theorem}\label{ammfiber}
Let $(R,\m)$ be a Cohen-Macaulay local ring of dimension $d \geq 2$.
Let $I$ be an $\m$-primary ideal with almost minimal multiplicity and
$\gamma(I) \geq d-2$. Then $F(I)$ has almost maximal depth.
\end{theorem}
\begin{proof}
By passing to $R[X]_{\m[X]}$, where $X$ is an indeterminate, we may
assume that $R/\m$ is infinite. Let $J$ be a minimal reduction of $I$.
We apply induction on $d$. Let $d = 2$. Suppose that $I$ has almost
minimal multiplicity. Then by Lemma \ref{genlemma2}, $\lambda(\m
I^n/\m JI^{n-1}) \leq 1$ for all $n \geq 1$.  Choose $x \in J$ such
that $x^o$ is superficial in $F(I)$ and $x^*$ is superficial in
$G(I)$. Let $`` - "$ denote image modulo $(x)$. Pick $y \in J$ such 
that $J = (x,y)$.
\vskip 2mm
\noindent
Then $\lambda(\bar{\m }\bar{I}/\bar{\m }\bar{J}) = 1$ and hence
$\lambda({\bar{\m }}\bar{I}^n/\bar{\m }\bar{J}\bar{I}^{n-1}) = 1 =
\lambda(\m I^n/\m JI^{n-1})$ for all $n = 1, \ldots, s$, where $s =
r_{\bar{J}}^{\bar{\m }}(\bar{I})$. 
For $j 
\geq 0$, consider the following exact sequence
$$
0 \longrightarrow \frac{\m I^j : x}{\m I^j : J}
\stackrel{y}{\longrightarrow} \frac{\m I^{j+1} : x}{\m I^j}
\stackrel{x}{\longrightarrow} \frac{\m I^{j+1}}{\m JI^j}
\longrightarrow \frac{\bar{\m }\bar{I}^{j+1}}{\bar{\m }\bar{J}\bar{I}^j}
\longrightarrow 0.
$$
\vskip 1mm
\noindent
{\it Claim :} $\m I^{j+1} : x = \m I^j$ for all $j \geq 0$.
\vskip 2mm
\noindent
For $j = 0$, it is clear. Suppose that $j > 0$ and that the claim is
true for $j-1$. Then $\m I^{j-1} \subseteq \m I^j : J \subseteq \m I^j
: x = \m I^{j-1}$, where the last equality follows by induction. Thus
$\m I^j : x = \m I^j : J$. Since, by Lemma \ref{redinvariance},
$r_J^\m (I) = s$, it follows from the the above exact sequence that
$\m I^{j+1} : x = \m I^j$ for all $j \geq 0$.  Therefore $x^o$ is a
regular element in $F(I)$ and $\depth F(I) \geq 1$.
\vskip 2mm
\noindent
Let $d > 2$. Choose $x_1, \ldots, x_d$ such that $J = (x_1, \ldots,
x_d)$, $x_1^*, \ldots, x_{d-2}^*$ is a regular sequence in $G(I)$ and
$x_1^o, \ldots, x_d^o$ is a superficial sequence in $F(I)$. Let $``-"$
denote images modulo $(x_1, \ldots, x_{d-2}).$ Then $\bar{R}$ is a
2-dimensional Cohen-Macaulay local ring. Since $x_1, \ldots, x_{d-2}$
is a part of minimal generating set for $I$, 
$$
\mu(\bar{I}) = \mu(I) - d + 2 = (e_0(I) - \lambda(R/I) + d - 1) - d +
2 = e_0(I) - \lambda(R/I) + 1 = e_0(\bar{I}) -
\lambda(\bar{R}/\bar{I}) + 1.
$$
Therefore $\bar{I}$ is
an ideal with almost minimal multiplicity in $\bar{R}$. Since $\dim 
\bar{R} = 2$, by the first part, $\depth F(\bar{I}) \geq 1$. Since
$(x_1^*, \ldots, x_{d-2}^*)$ is a regular sequence in $G(I)$,
$F(\bar{I}) \cong F(I)/(x_1^o, \ldots, x_{d-2}^o)F(I)$
and hence by Proposition \ref{supeltsprprties}(4), $\depth F(I)
\geq d-1$.
\end{proof}
\vskip 2mm
\noindent
We end this section with an example to show that the depth assumption
on the associated graded ring in the above theorem is necessary. This 
example was provided to us by M. E. Rossi.
\begin{example}{\em 
Let $R = k[\![x, y, z]\!]$, where $k$ is any field.  Let $I =
(-x^2+y^2,-y^2+z^2, xy, yz, zx)$ and $J = (-x^2+y^2, -y^2+z^2, xy)$.
Then $I^3 = JI^2$. Hence $J$ is a minimal reduction of $I$. Let $\m =
(x,y,z)$. Then it can be seen that $\m I = \m J + (z^3)$ and $\m(z^3)
\subset \m J$. Hence $\lambda(\m I/\m J) = 1$. Therefore
$I$ has almost minimal multiplicity. It can be easily seen that $x^2I
\subset I^2$, but $x^2 \notin I$. This shows that the Ratliff-Rush
closure $\tilde{I}$ is not equal to $I$. Hence, $\depth G(I) = 0$.

Now we show that $\depth F(I) = 1$.  Since I is generated by
homogeneous elements of same degree (equal to 2), $F(I) \cong
k[-x^2+y^2,-y^2+z^2, xy, yz, zx]$. Therefore $\depth F(I) \geq 1$. The
minimal generating set, $(-x^2+y^2, -y^2+z^2, xy)$ of $J$ form a
system of parameters for $F(I)$. If $\depth F(I) > 1$, then there
exist at least one pair in this set such that they form a regular
sequence in $F(I)$. Let $F = k[-x^2+y^2,-y^2+z^2, xy, yz, zx]$. Then 
$(y^2z^2 - x^2z^2)(-y^2+z^2) = (-x^2+y^2)(-y^2z^2 + z^4)$. Since
$-y^2z^2+z^4 \in F$, $(-x^2+y^2)F : (-y^2+z^2)F = (-x^2+y^2,
-x^2z^2+y^2z^2)F$. Also, note that $y^2z^2 - x^2z^2 \notin
(-x^2+y^2)F$. This shows that $-y^2+z^2$ is a zerodivisor in $F$
$\modulo(-x^2+y^2)F$. Similarly, one can see that for any
pair of elements, $t_1, t_2 \in \{-x^2+y^2, -y^2+z^2, xy\}, \; \;
t_1F : t_2F \neq t_1F$. Hence $\depth F = \depth F(I) = 1$. 

}
\end{example}
\section{Cohen-Macaulay $F_K(I)$ when $I$ has almost minimal
multiplicity}
In this section, we characterize Cohen-Macaulay property of $F_K(I)$
when $I$ has almost minimal multiplicity. For this purpose, we find
the generating function of the function $H_K(I, n)$, first in
dimension 1 and then in arbitrary dimension. A formula of Rossi and
Valla for the Hilbert series of $G(\m)$ when $\lambda(\m^2/J\m) = 1$ is
generalized for $\m$-primary ideals with almost minimal multiplicity. 

\begin{lemma}\label{polyseries}
Let $(R,\m)$ be a $1$-dimensional Cohen-Macaulay local ring and let
$I$ be an $\m$-primary ideal of $R$ with almost minimal multiplicity
with respect to $K \supseteq I$. Let $s = r_J^K(I)$. Then 
\begin{enumerate}
\item[1.] $P_{K}(I,n) = e_0(I)n - (s-\lambda(R/K))$.
\item[2.] $\sum_{n\geq0}H_{K}(I,n)t^n = \left[\lambda(R/K) + 
(e_0(I)-1-\lambda(R/K))t + t^{s+1}\right]/(1-t)^{2}.$ 
\end{enumerate}
\end{lemma}
\begin{proof} 
1. From the following diagram, 
$$
\xymatrix{
 R \ar@{<-}[r]^h \ar@{<-}[d]_{e_0(I)}
   & KI_n \ar@{<-}[r]^{l_n} \ar@{<-}[d]_{e_0(I)} & KI_{n+1}
   \ar@{<-}[dl]^{c} \\
    (x) \ar@{<-}[r]_h & xKI_n &           }
$$
it follows that 
\begin{eqnarray*}
\lambda(KI^n/KI^{n+1}) & = & \left\{
                               \begin{array}{ll}
			          e_0(I) - 1 \hspace*{0.1in}\mbox{ for
				  } n = 0, \ldots, s-1 \\
				  e_0(I) \hspace*{0.1in}\mbox{ for } n
				  \geq s.
				\end{array}\right.
\end{eqnarray*}

Therefore, 
\begin{eqnarray*}
\lambda(R/KI^n) & = & \lambda(R/K) +
\sum_{i=0}^{n-1}\lambda(KI^n/KI^{n+1}) \\
 & = & \left\{
\begin{array}{ll}
n(e_0(I)-1) + \lambda(R/K) \mbox{ for } 1 \leq n \leq s \\
ne_0(I) - (s - \lambda(R/K)) \mbox{ for } n > s.
\end{array}\right.
\end{eqnarray*}
Therefore $P_{K}(I,n) = ne_0(I) - (s- \lambda(R/K))$.

\vskip 4mm
\noindent
2. Substituting the values of $H_{K}(I,n)$ from (1) we get, 
\begin{eqnarray*}
\sum_{n\geq0} H_{K}(I,n)t^n & = & \sum_{n=0}^s
[n(e_0(I)-1)+\lambda(R/K)]t^n + \sum_{n=s+1}^{\infty} [e_0(I)n -
(s-\lambda(R/K))]t^n \\
& = & e_0(I) \sum_{n=0}^{\infty} nt^n - \sum_{n=0}^s nt^n + 
\lambda(R/K) \sum_{n=0}^{\infty} t^n - s \sum_{n=s+1}^{\infty}t^n\\
& = & \frac{e_0(I)}{(1-t)^2} - \frac{e_0(I)}{(1-t)} + 
\frac{\lambda(R/K)-st^{s+1}}{(1-t)} - \sum_{n=0}^s nt^n\\
& = & \frac{e_0(I)}{(1-t)^2} - \frac{e_0(I) - \lambda(R/K) + st^{s+1}+
(\sum_{n=0}^s nt^n)(1-t)}{(1-t)} \\
& = & \frac{e_0(I)}{(1-t)^2} -
\frac{e_0(I) - \lambda(R/K) + st^{s+1} + t(1+t+ \cdots +t^{s-1}) -
st^{s+1}}{(1-t)}\\
& = & \frac{e_0(I) - e_0(I)(1-t) + \lambda(R/K)(1-t) - t(1-t^s)}
{(1-t)^2} \\
& = & \frac{\lambda(R/K)+ (e_0(I)-1-\lambda(R/K))t + t^{s+1}}
{(1-t)^2}.
\end{eqnarray*}
\end{proof}
\begin{proposition}\label{hilbser}
Let $(R,\m)$ be a Cohen-Macaulay local ring of dimension $d$. Let $I$
be an $\m$-primary ideal with almost minimal multiplicity with respect
to $K$ such that $\gamma(I) \geq d-1$. Let $s = r_J^{K}(I)$. Then
$$
\sum_{n\geq0}H_{K}(I,n)t^n = \frac{\lambda(R/K) + (e_0(I) - 1- 
\lambda(R/K))t + t^{s+1}}{(1-t)^{d+1}}.
$$
\end{proposition}
\begin{proof} We induct on $d$. The case $d = 1$ is proved in Lemma
\ref{polyseries}(2). 
Let $d > 1$.  Let $x \in I \backslash KI$,
such that $x^*$ is a regular element in $G(I)$ and $x^o$ is a regular
element in $F_K(I)$. Let $``-"$ denote images modulo $(x)$.  Then
$\bar{I}$ is an $\bar{\m}$-primary ideal with almost minimal
multiplicity with respect to $\bar{K}$ in $\bar{R}$.  For $n \geq 1$,
consider the exact sequence
$$
\begin{CD}
0 @>>> KI^{n+1} : x/KI^n @>>> R/KI^n @>x>> R/KI^{n+1} @>>>
R/(KI^{n+1}+xR) @>>> 0.
\end{CD}
$$
Since $x^*$ is regular in $G(I)$ and $x^o$ is regular in $F_K(I)$, 
$KI^{n+1} : x = KI^n$ for all $n \geq 0$, by Proposition
\ref{supeltsprprties}(3).
Therefore $H_{\bar{K}}(\bar{I}, n) = \Delta H_{K}(I,n)$ for all $n 
\geq 1$. By induction 
$$
\sum_{n\geq0} H_{\bar{K}}(\bar{I}, n)t^n =
\frac{\lambda(\bar{R}/\bar{K}) + (e_0(\bar{I}) - 1 - 
\lambda(\bar{R}/\bar{K}))t + t^{s+1}}{(1-t)^d},
$$
where $s = r_{\bar{J}}^{\bar{K}}(\bar{I}) = r_J^{K}(I)$. 
Therefore 
$$
\sum_{n\geq0}H_{K}(I,n)t^n = \frac{\lambda(R/K) + (e_0(I) - 1 -
\lambda(R/K))t + t^{s+1}}{(1-t)^{d+1}}.
$$
\end{proof}
As a consequence of the above result, we recover a result of M. E.
Rossi, \cite[Corollary 3.8(2)]{r1}

\begin{corollary}
Let $(R,\m)$ be a Cohen-Macaulay local ring of dimension $d > 0$. Let
$I$ be an $\m$-primary ideal with $\lambda(I^2/JI) = 1$ for some
minimal reduction $J$ of $I$ with reduction number $r$. Then 
$$
H(G(I),t) := \sum_{n\geq0}\lambda(I^n/I^{n+1})t^n =
\frac{\lambda(R/I) + (e_0(I) - 1 -\lambda(R/I))t + t^r}{(1-t)^d}.
$$
\end{corollary}
\begin{proof} By Corollary 1.7 of \cite{r}, $\depth G(I) \geq d-1$.
Put $K = I$ in Proposition \ref{hilbser}. Then we get 
$$
\sum_{n\geq0}\lambda(R/I^{n+1})t^n = \frac{\lambda(R/I) + (e_0(I) - 1 - 
\lambda(R/I))t + t^r}{(1-t)^{d+1}}.
$$
Multiplying both sides by $(1-t)$, we get
$$
\sum_{n\geq0}\lambda(I^n/I^{n+1})t^n = \frac{\lambda(R/I) + (e_0(I) -
1 -\lambda(R/I))t + t^r}{(1-t)^d}.
$$
\end{proof}

\vskip 2mm
\noindent
We end this paper by characterizing the Cohen-Macaulay fiber cones of
ideals with almost minimal multiplicity in the following proposition. 

\begin{proposition}\label{ammcmchar}
Let $(R,\m)$ be a $d$-dimensional Cohen-Macaulay local ring and $I$ be
an $\m$-primary ideal with almost minimal multiplicity with respect to
$K \supseteq I$ and $\gamma(I) \geq d-1$. Let $s = r_J^{K}(I)$. Then
$F_K(I)$ is Cohen-Macaulay if and only if $\lambda(KI^n +
JI^{n-1}/JI^{n-1}) = 1$ for all $n = 1, \ldots, s$.
\end{proposition}
\begin{proof}
Since $I$ has almost minimal multiplicity with respect to $K$ and 
$\gamma(I) \geq d-1$, it follows from the generating function given in
Proposition \ref{hilbser} and Proposition 4.1.9 of \cite{bh} that
$g_1 = s-\lambda(R/K)$. As in the proof of Lemma \ref{genlemma2}, it
can easily be seen that $\lambda(KI^n+JI^{n-1}/JI^{n-1}) \leq
1$ for all $n \geq 1$. Therefore $F_K(I)$ is CM if and only
if
$$
\sum_{n=1}^{s}\lambda(KI^n+JI^{n-1}/JI^{n-1}) - \lambda(R/K) = s -
\lambda(R/K).
$$
if and only if $\lambda(KI^n+JI^{n-1}/JI^{n-1}) = 1$ for all $n = 1,
\ldots, s.$
\end{proof}

\end{document}